\documentstyle[amssymb,amsfonts]{amsart}

\def\R{{\hbox{\bf R}}}

\def\Q{{\hbox{\bf Q}}}

\def\X{{E}}

\font \roman = cmr10 at 10 true pt

\def\p{{\hbox{\bf p}}}
\def\E{{\cal E}}
\def\A{{\cal A}}

\def\allt#1{%
\smash{
\vtop{%
     \ialign{%
        ##\crcr
        $\hfil\displaystyle{\tilde \forall}\hfil$\crcr%
        \noalign{\kern1.5pt\nointerlineskip}
        $\hfil\!\!#1\hfil$\crcr\noalign{\kern1.5pt}
        }
       }
      } \hbox{$\vphantom{#1}$}
     }

\def\be#1{\begin{equation}\label{#1}}
\def\bas{\begin{align*}}
\def\eas{\end{align*}}
\def\bi{\begin{itemize}}
\def\ei{\end{itemize}}

\def\dist{{\hbox{\roman dist}}}
\def\rp{{r^\prime}}

\def\qp{{q^\prime}}

\def\eps{\varepsilon}
\newenvironment{proof}{\noindent {\bf Proof} }{\endprf\par}
\def \endprf{\hfill  {\vrule height6pt width6pt depth0pt}\medskip}
\def\emph#1{{\it #1}}
\def\textbf#1{{\bf #1}}

\def\calg{{\cal G}}

\theoremstyle{plain}
  \newtheorem{theorem}[subsection]{Theorem}

  \newtheorem{proposition}[subsection]{Proposition}
  
  \newtheorem{lemma}[subsection]{Lemma}
  \newtheorem{corollary}[subsection]{Corollary}

\theoremstyle{remark}

\theoremstyle{definition}

\include{psfig}

\begin{document}

\title[An x-ray transform estimate in $\R^n$]{An x-ray transform estimate in $\R^n$}

\author{Izabella {\L}aba}
\address{Department of Mathematics, Princeton University, Princeton NJ
08544}
\email{laba@@math.princeton.edu}

\author{Terence Tao}
\address{Department of Mathematics, UCLA, Los Angeles CA 90095-1555}
\email{tao@@math.ucla.edu}


\begin{abstract}  We prove an x-ray estimate in general dimension which is a stronger version of Wolff's Kakeya estimate \cite{W1}.  This generalizes the estimate in \cite{W2}, which dealt with the $n=3$ case.
\end{abstract}

\maketitle

\section{Introduction}

Let $n \geq 3$ be an integer.  Let $B^{n-1}(0,1)$ be the unit ball in $\R^n$, and for all $x, v \in B^{n-1}(0,1)$ define the line segment $l(x,v) \in \R^n$ by
$$ l(x,v) = \{ (x+vt, t): t \in [0,1] \}$$
where we have parameterized $\R^n$ as $\R^{n-1} \times \R$ in the usual manner.  
Let $\calg$ be the set of all such line segments; this space is thus identified with $B^{n-1}(0,1) \times B^{n-1}(0,1)$.  If $l \in \calg$, we write $x(l)$ and
$v(l)$ for the values of $x$ and $v$ respectively such that $l = l(x,v)$.

For any function $f$ on $\R^n$, define the \emph{x-ray transform} $Xf$ on $\calg$ by
$$ X f(l) = \int_l f.$$
We consider the question of determining the exponents $1 \leq p,q,r \leq \infty$ and $\alpha \geq 0$ such that we have the bound
\be{x-ray}
 \| Xf \|_{L^q_v L^r_x} \lesssim \| f\|_{L^p_\alpha},
\end{equation}
where $L^p_\alpha$ is the Sobolev space $(1 + \sqrt{-\Delta})^{-\alpha} L^p$.

From scaling considerations (or by letting $f$ be a bump function adapted to a small ball) we have the necessary condition
\be{scaling}
1 + \frac{n-1}{r} \geq \frac{n}{p} - \alpha
\end{equation}
while if one lets $f$ be adapted to a tubular neighbourhood of a line segment $l \in \calg$, we obtain the condition
\be{knapp}
\frac{n-1}{q} + \frac{n-1}{r} \geq \frac{n-1}{p} - \alpha.
\end{equation}
From the Besicovitch set construction we have
\be{besicovitch}
(r,\alpha) \neq (\infty, 0).
\end{equation}

It was conjectured by Drury \cite{drury:x-ray} and Christ \cite{christ:x-ray} that these three necessary conditions are in fact sufficient. In \cite{christ:x-ray} this conjecture was shown to be true when $p \leq (n+1)/2$.

By H\"older, Sobolev, and interpolation with trivial estimates, the full conjecture is equivalent (modulo endpoints) to the Kakeya conjecture, which asserts that \eqref{x-ray} holds for $q = n$, $r = \infty$, $p = n$, and $\alpha = \eps$ for arbitrarily small $\eps$.

Wolff \cite{W1} showed \eqref{x-ray} was true when 
$$q = \frac{(n-1)(n+2)}{n}, \quad r = \infty, \quad p = \frac{n+2}{2}, \quad \alpha = \frac{n-2}{n+2} + \eps;$$
this can of course be interpolated with the results in \cite{christ:x-ray} to yield further estimates.  However, this is not the best one can do in the $p = \frac{n+2}{2}$ case.  From \eqref{scaling} and \eqref{knapp} one expects to have \eqref{x-ray} for
\be{endpt}
q = \frac{(n-1)(n+2)}{n}, \quad r = \frac{(n-1)(n+2)}{n-2}, \quad p = \frac{n+2}{2}, \quad \alpha = 0;
\end{equation}
this would imply the results of \cite{W1} by Sobolev embedding in the $v$ variable.  Although we are not able to get that sharp result, we are able to
obtain the following interpolant, which is our main result.

\begin{theorem}\label{main}  For any $\eps > 0$, we have \eqref{x-ray} for
\be{squid}
q = \frac{(n-1)(n+2)}{n}, \quad r = 2(n+2), \quad p = \frac{n+2}{2}, \quad \alpha = \frac{n-3}{2(n+2)} + \eps.
\end{equation}
\end{theorem}

This result was obtained in the three dimensional case $n=3$ by Wolff \cite{W2},
and the result is sharp up to endpoints for that value of $n$ and $p$.  Our arguments shall be based on those in \cite{W2}, with some mild simplifications based on the bilinear approach in \cite{TVV}.

Theorem \ref{main} can be stated in a discretized adjoint form, which is more convenient for applications.  Namely\footnote{The notation in the theorem will be explained shortly.}:

\begin{theorem}\label{discrete-thm}  Let $\eps > 0$, $0 < \delta \ll 1$, and
$1 \leq m \lesssim \delta^{1-n}$.
Let $\E$, $\E'$ be $\delta$-separated subsets of $B^{n-1}(0,1)$, and let
$\A \subset \E \times \E' \subset \calg$ be a collection of line segments such that
\be{m-def}
|\{ l \in \A: v(l) = v \}| \leq m
\end{equation}
for all $v \in \E$.  Then we have
\be{discrete-est}
\| \sum_{l \in \A} \chi_{T_l} \|_{p'} \lesssim 
\delta^{-\frac{n}{p} + 1 - \eps} 
m^{1/q - 1/r} (\delta^{n-1} |\A|)^{1/\qp}
\end{equation}
where $p,q,r$ are as in \eqref{squid}.
\end{theorem}

As observed in \cite{W2}, an x-ray estimate of this form reveals some information on Besicovitch sets in $\R^n$.  Namely, such sets have Minkowski dimension at least $\frac{n+2}{2}$, and if the dimension is exactly $\frac{n+2}{2}$ then the line segments which comprise the set must be ``sticky''
in a certain sense.  This observation was made rigorous in \cite{KLT}, where 
the results of \cite{W2} were applied (together with those of \cite{B3} and some additional arguments) in the three-dimensional case to improve slightly upon the Minkowski bound just stated.  We will use Theorem \ref{main} to achieve a similar result in higher dimensions \cite{LT}.  Fortunately, one does not need a sharp value of $r$ in \eqref{squid} to obtain this type of observation, as long as $r$ is finite of course.

To illustrate the connection between x-ray estimates and Besicovitch sets, we note the following simple application of Theorem \ref{discrete-thm}:

\begin{corollary}  Let $0 \leq \alpha \leq n-1$, and 
$E$ be a bounded subset of $\R^n$ such that, for each direction $\omega \in S^{n-1}$, $E$ contains a family of unit line segments parallel to $\omega$,
whose union has Minkowski dimension $\alpha + 1$.  Then the Minkowski dimension of $E$ is at least $\frac{n+2}{2} + \frac{\alpha}{4}$.
\end{corollary}

The proof follows standard discretization arguments (see e.g. \cite{B1}, \cite{B2}) and will
be omitted.  A similar result holds when Minkowski dimension is replaced by Hausdorff.  This corollary is stronger than the corresponding corollary of the Kakeya estimate in \cite{W1}, which covers the $\alpha = 0$ case.  If one had an x-ray estimate for \eqref{endpt} then one would be able to improve the $\alpha/4$ term to the optimal $\frac{n-2}{2n-2} \alpha$.

The second author is supported by NSF grant DMS-9706764, and wishes to thank Nets Katz and Tom Wolff for helpful discussions.

\section{Notation}

We use $0 < \delta \ll 1$ and $0 < \eps \ll 1$ to denote certain small numbers, and $N \gg 1$ denotes a certain large integer.  If $l$ is a line segment in $\calg$, we use $T_l$ to denote the $\delta$-neighbourhood of $l$, which is thus a $\delta \times 1$ tube.

We write $A \lesssim B$ for $A \leq CB$, $A \ll B$ for $A \leq C^{-1} B$, and $A \lessapprox B$ for
$A \leq C (\log(1/\delta))^\nu B$, and $C, \nu$ are quantities which vary from line to line and are allowed to depend on $\eps$ and $N$ but not on $\delta$.  

Our argument will require the introduction of many quantities, which measure various angles or cardinalities in a collection of tubes.  For purposes of visualizing the argument we recommend that one sets the values of these quantities as follows:
$$ |\E| \sim |\E_i| \sim \delta^{1-n},\ |\A| \sim \delta^{1-n} m,\ \lambda \sim \theta \sim \sigma \sim 1,\ \rho \sim \p_i \sim w$$
for $i=1,2$.  The treatment of this case can be done while avoiding the more technical tools in the argument such as the two-ends and bilinear reductions, and most of the uniformization theory, while still capturing the core ideas of
the argument.  To improve the value of $r$ in \eqref{squid} one would probably start by considering this case.

\section{Derivation of Theorem \ref{main} from Theorem \ref{discrete-thm}}\label{discrete}

Assume that Theorem \ref{discrete-thm} holds.  In this section we shall see how Theorem \ref{main} follows.  The argument is standard (cf. \cite{B1}, \cite{B2}, \cite{W1}, \cite{W2}, \cite{TVV}).

By a Littlewood-Paley decomposition, and giving up an epsilon in the $\alpha$ index, one may assume that $f$ has Fourier support in an annulus $\{ \xi: |\xi| \sim \delta^{-1} \}$.  The case $\delta \gtrsim 1$ is easy to handle, so we assume henceforth that $0 < \delta \ll 1$.

Fix $\delta$.  It is then well known that \eqref{x-ray} follows from the variant
$$ 
\| X_\delta f \|_{L^q_v L^r_x} \lessapprox \delta^{-\alpha} \| f\|_{p},
$$
where
$$ X_\delta f(l) = \delta^{1-n} \int_{T_l} f$$
By duality this is equivalent to
$$ \| X_\delta^* F \|_{p'} \lessapprox \delta^{-\alpha} \| F \|_{L^\qp_v L^\rp_x}$$
for all $F$ on $\calg$, where $X_\delta^*$ is the adjoint x-ray transform
$$ X_\delta^* F = \delta^{1-n} \int_\calg F(l) \chi_{T_l}\ dx dv.$$

Let $\E$, $\E'$ by any $\delta$-separated subsets of $B^{n-1}(0,1)$.  By discretization it suffices to show that
$$ \| \delta^{n-1} \sum_{v \in \E} \sum_{x \in \E'} F(l(x,v)) \chi_{T_{l(x,v)}} 
\|_{p'} \lessapprox \delta^{-\alpha} ( \delta^{n-1} \sum_{v \in \E}
( \delta^{n-1} \sum_{x \in \E'} |F(l(x,v))|^\rp)^{\qp/\rp} )^{1/\qp}$$
uniformly in $\E$, $\E'$.

Fix $\E$, $\E'$
By pigeonholing and positivity it suffices to verify this when $F$ is a characteristic function $F = \chi_\A$ for some $A \subseteq \E \times \E'$, so that we reduce to
$$ \| \sum_{l \in \A: v(l) \in \E} \chi_{T_l} 
\|_{p'} \lessapprox \delta^{(n-1)(1 - \frac{1}{r} - \frac{1}{q}) - \alpha}
( \sum_{v \in \E} |\{ l \in \A: v(l) = v \}|^{\qp/\rp} )^{1/\qp}.$$
By a further pigeonholing and refining of $\E$, we may assume that there exists $1 \leq m \lessapprox \delta^{1-n}$ such that 
\be{m-def-ex}
m/2 \leq |\{ l \in \A: v(l) = v \}| \leq m
\end{equation}
for all $v \in \E$.
Our task is then to show that
$$ \| \sum_{l \in \A: v(l) \in \E} \chi_{T_l} 
\|_{p'} \lessapprox \delta^{(n-1)(1 - \frac{1}{r} - \frac{1}{q}) - \alpha} m^{1/\rp} |\E|^{1/\qp}.
$$
From \eqref{m-def-ex} we then have $|\A| \sim m |\E|$.  The claim then follows from
Theorem \ref{discrete-thm} and the fact that \eqref{scaling} is almost satisfied with equality.

It thus remains to prove Theorem \ref{discrete-thm}.

\section{A three-dimensional estimate}

For any collection $\A$ of line segments, we follow Wolff \cite{W2} (see also \cite{W3}) and define the \emph{plate number} $\p(\A)$ by
\be{plate-def}
 \p(\A) = \sup_R \frac{|\{ l \in \A: T_l \in R\}|}{w/\delta}
\end{equation}
where $R$ ranges over all rectangles of dimension $C \times Cw \times C\delta \times \ldots \times C\delta$.  By considering the $w \sim \delta$ case we see that $\p(\A) \gtrsim 1$ for any non-empty $\A$.

The purpose of this section is to prove the following distributional estimate on a set $\X$ assuming that the directions of $\A$ are effectively constrained to a two-dimensional slab, and the intersection of the tubes $T_l$ with $\X$ satisfy a certain ``two-ends'' condition of the type used in \cite{W1}, \cite{W2}.  This lemma will be key in the main argument, and also
employs several techniques, notably a hairbrush argument
and a uniformization argument (both due to Wolff), which will re-appear in
slightly different form in the sequel.

\begin{lemma}\label{three}\cite{W2}
Let $N \gg 1$ be an integer, $\delta^C \lesssim \lambda \leq 1$, $\X$ be a subset of $\R^n$, and let $\A \subset \E \times \E'$ be a collection of lines satisfying \eqref{m-def} which satisfy the uniform density estimate
\be{lam-as} |T_l \cap \X| \approx \lambda \delta^{n-1}
\end{equation}
and the two-ends condition
\be{two-ends-three}
 |T_l \cap \X \cap B(x, \delta^{1/N})| \lessapprox \delta^{\eps/2N}
\lambda \delta^{n-1}
\end{equation}
for all $l \in \A, x \in \R^n$.  Suppose also that the set of directions
$\{ v(l): l \in \A\}$ is contained in a $C \times C\rho \times C\delta \times \ldots \times C\delta$ box in $B^{n-1}(0,1)$ for some $\delta \lesssim \rho \lesssim 1$.  Then, if $\delta$ is sufficiently small depending on $\eps$ and $N$, we have
\be{xa-targ} |\X| \gtrapprox \delta^{C/N}
\lambda^{2} |\A|
m^{-1/2} \rho^{-1/2} \p(\A)^{-1/2} \delta^{n-1/2}.
\end{equation}
\end{lemma}

\begin{proof}
We repeat the argument in \cite{W2}.  We may assume that $\A$ is non-empty, and that $\X$ is contained in $\bigcup_{l \in \A} T_l$.

For every $l \in \A$ and dyadic $\delta \lesssim \sigma \lesssim 1$, $1 \leq \mu \lesssim \delta^{-C}$, we let
$Y_{l,\mu,\sigma,\A} \subset T_l \cap \X$ denote the set
\be{y-def}
 Y_{l,\mu,\sigma,\A} = \{ x \in T_l \cap \X: \sum_{l' \in \A: \delta + |v(l)-v(l')| \sim \sigma} \chi_{T_{l'}}(x) \approx  \sum_{l' \in \A} \chi_{T_{l'}}(x) \approx \mu\}.
\end{equation}
In other words, $Y_{l,\mu,\sigma,\A}$ consists of those points $x$ in $T_l \cap X$ which lies in about $\mu$ tubes from $\A$, most of which make an angle of about $\sigma$ with $T_l$.  From the pigeonhole principle we see that
\be{wav}
 T_l \cap \X = \bigcup_{\delta \lesssim \sigma \lesssim 1} \bigcup_{1 \leq \mu \lesssim \delta^{-C}} Y_{l,\mu,\sigma,\A}.
\end{equation}

We now prove a technical lemma which allows us to uniformize $\mu$ and $\sigma$.  This type of argument will also be used in the sequel.  (For a more general formulation of this type of argument, see \cite{W2}).  A somewhat similar lemma appears in \cite{ccc}.

\begin{lemma}\label{tech}  Let the notation be as above.  Then there exist quantities
$\delta \lesssim \sigma \lesssim 1$, $1 \leq \mu \lesssim \delta^{-C}$
and sets
$$ \A^{(2)} \subseteq \A^{(1)} \subseteq \A^{(0)} = \A$$
and for each $i = 1,2$, $l \in \A^{(i)}$ there exists a set
$$ Y_l^{(i)} \subseteq T_l \cap \X$$
such that
\be{ai-size}
 |\A^{(i)}| \approx |\A|,
\end{equation}
\be{yli-size}
 |Y_l^{(i)}| \approx \lambda \delta^{n-1},
\end{equation}
and 
\be{yi-unif}
 Y_l^{(i)} \subseteq Y_{l,\mu^{(i)},\sigma^{(i)},\A^{i-1}}
\end{equation}
for some $\mu^{(i)}, \sigma^{(i)}$ satisfying
\be{mu-unif-3}
 \delta^{C/N} \mu \lesssim \mu^{(i)} \lesssim \delta^{-C/N} \mu
\end{equation}
\be{sigma-unif-3}
 \delta^{C/N} \sigma \lesssim \sigma^{(i)} \lesssim \delta^{-C/N} \sigma
\end{equation}
The implicit constants may depend on $N$.
\end{lemma}

\begin{proof}  The first stage shall be to construct sequences
$$ \A = \A_0 \supseteq \A_1 \supseteq \ldots \supseteq \A_{N^2},$$
$$ T_l \cap \X = Y_{l,0} \supset Y_{l,1} \supseteq \ldots \supseteq Y_{l,k},$$
and quantities $\mu_k, \sigma_k$ for all $1 \leq k \leq N^2$ and $l \in \A_k$, such that
\be{a-card-iter}
|\A_k| \approx |\A|,
\end{equation}
\be{y-size-iter}
|Y_{l,k}| \approx \lambda \delta^{n-1},
\end{equation}
and
\be{y-angle-iter}
Y_{l,k} \subseteq Y_{l,\mu_k,\sigma_k,\A_{k-1}}
\end{equation}
for all $1 \leq k \leq N^2$.

To do this, suppose inductively that $0 \leq k < N^2$ is such that $\A_k$ and
$Y_{l,k}$ have been constructed for all $l \in \A_k$.  From \eqref{wav} we have
$$
Y_{l,k} \subseteq \bigcup_{\delta \lesssim \sigma \lesssim 1} \bigcup_{1 \leq \mu \lesssim \delta^{-C}} Y_{l,\mu,\sigma,\A_k}.$$
By the pigeonhole principle, for every $l \in \A_k$ one can thus find
$\mu_k(l)$, $\sigma_k(l)$ such that
$$ |Y_{l,k+1}| \approx |Y_{l,k}|,$$
where
$$ Y_{l,k+1} = Y_{l,k} \cap Y_{l,\mu_k(l),\sigma_k(l),\A_k}.$$
By the pigeonhole principle again, there exists $\mu_k$, $\sigma_k$ independent of $l$ such that the set
$$ \A_{k+1} = \{ l \in \A_k: \mu_k(l) = \mu_k, \sigma_k(l) = \sigma_k \}$$
satisfies \eqref{a-card-iter}.  It is clear that this construction gives the desired properties.

By the pigeonhole principle, there must exist $1 \leq k_1 < k_2 \leq N^2$
and $\sigma$, $\mu$ such that
$$ \delta^{C/N} \mu \lessapprox \mu_{k_i} \lessapprox \delta^{-C/N} \mu$$
and
$$ \delta^{C/N} \sigma \lessapprox \sigma_{k_i} \lessapprox \delta^{-C/N} \sigma$$
for $i=1,2$.  The claim then follows by setting $\A^{(i)} = \A_{k_i}$
and $Y_l^{(i)} = Y_{l,k_i}$.
\end{proof}

Let the notation be as in the above lemma.  From \eqref{yli-size} and \eqref{ai-size} we have
$$ \sum_{l \in \A^{(2)}} |Y_l^{(2)}| \approx \lambda \delta^{n-1} |\A|$$
which we rewrite as
$$ \int_\X \sum_{l \in \A^{(2)}} \chi_{Y_l^{(2)}} \approx \lambda \delta^{n-1} |\A|.$$
From \eqref{yi-unif}, the nesting $\A^{(2)} \subseteq \A^{(1)}$, and \eqref{mu-unif-3}, the integrand is bounded by $\delta^{-C/N} \mu$.  We thus see that $\lambda$ and $\mu$ are naturally related by the estimate
\be{mu-lambda-3}
|\X| \mu \gtrapprox \delta^{C/N} \lambda \delta^{n-1} |\A|.
\end{equation}
One can reverse the inequality in \eqref{mu-lambda-3}, but we shall not need to do so here.

From \eqref{a-card-iter}, $\A^{(2)}$ is non-empty.
Let $l_0$ be an arbitrary element of $\A^{(2)}$.  Consider the ``hairbrush'' $\A_{brush}^{l_0}$ defined by
$$ \A_{brush}^{l_0} = \{ l \in \A^{(1)}: T_{l_0} \cap T_{l} \neq \emptyset,
\delta^{C/N} \sigma \lessapprox
\delta + |v(l_0)-v(l)| \lessapprox \delta^{-C/N} \sigma \}.$$
From \eqref{yi-unif}, \eqref{mu-unif-3}, \eqref{sigma-unif-3} we see that 
$$ \sum_{l \in \A_{brush}^{l_0}} \chi_{T_l}(x) \gtrapprox \delta^{C/N} \mu$$
for all $x \in Y_{l_0}^{(2)}$.
Integrating this using \eqref{yli-size}, we obtain
$$ \sum_{l \in \A_{brush}^{l_0}} |T_l \cap Y_{l_0}^{(2)}| 
\gtrapprox \delta^{C/N} \mu \lambda \delta^{n-1}.$$
From elementary geometry we see that
$$ |T_l \cap Y_{l_0}^{(2)}| \leq |T_l \cap T_{l_0}| \lessapprox \delta^{-C/N} \delta^n \sigma^{-1}$$
so we conclude that
\be{abrush-card}
|\A_{brush}^{l_0}| \gtrapprox \delta^{C/N} \mu \lambda \sigma \delta^{-1}.
\end{equation}

We will shortly combine \eqref{abrush-card} with \eqref{yli-size} and
\eqref{two-ends-three} to prove the estimate
\be{brush-bound}
|\bigcup_{l \in \A_{brush}^{l_0}} Y_l^{(1)}| \gtrapprox
\delta^{C/N} \mu \lambda^3 \sigma \p(\A)^{-1} \delta^{n-2}.
\end{equation}

Assuming this bound for the moment, let us complete the proof of \eqref{xa-targ}.  From \eqref{yi-unif} and \eqref{sigma-unif-3} we have
$$ \sum_{l' \in \A: \delta + |v(l)-v(l')| \lesssim \delta^{-C/N} \sigma} \chi_{T_{l'} \cap \X}(x) \gtrapprox \delta^{C/N} \mu$$
for all $l \in \A_{brush}^{l_0}$ and $x \in Y_l^{(1)}$.  From the 
definition of $\A_{brush}^{l_0}$ and the triangle inequality 
we thus see that
$$ \sum_{l' \in \A: \delta + |v(l')-v(l_0)| \lessapprox \delta^{-C/N} \sigma} \chi_{T_{l'} \cap \X}(x) \gtrapprox \delta^{C/N} \mu$$
for all $x$ in the set in \eqref{brush-bound}.  Integrating this and using \eqref{brush-bound}, we thus obtain
$$ \sum_{l' \in \A: \delta + |v(l')-v(l_0)| \lessapprox \delta^{-C/N} \sigma} |T_{l'} \cap \X| \gtrapprox \delta^{C/N} \mu^2 \lambda^3 \sigma \p(\A)^{-1} \delta^{n-2}$$
From \eqref{lam-as} we thus have
$$ |\{ l' \in \A: \delta + |v(l')-v(l_0)| \lessapprox \delta^{-C/N} \sigma \}|
\lambda \delta^{n-1} \gtrapprox \delta^{C/N} \mu^2 \lambda^3 \sigma 
\p(\A)^{-1} \delta^{n-2}.$$
However, from \eqref{m-def} and the fact that $v(l')$ is constrained to a
$C \times C\rho \times C\delta \times \ldots \times C\delta$ box, we see from elementary geometry that
$$
|\{ l' \in \A: \delta + |v(l')-v(l_0)| \lessapprox \delta^{-C/N} \sigma \}|
\lessapprox \delta^{-C/N} \sigma \rho \delta^{-2}.$$
Combining these two estimates we obtain (after some algebra)
$$ \mu \lessapprox \delta^{-C/N} \rho^{1/2} \delta^{-1/2} \lambda^{-1}
\p(\A)^{1/2},$$
and the claim \eqref{xa-targ} follows after some algebra from this and
\eqref{mu-lambda-3}.

It remains to prove \eqref{brush-bound}.  We first deal with a trivial case when $\sigma \lesssim \delta^{-C/N} \delta$.  In this case we simply use the bound
$$ |\bigcup_{l \in \A_{brush}^{l_0}} Y_l^{(1)}| \geq Y_l^{(1)}
\gtrsim \delta^{C/N} \lambda \delta^{n-1}$$
from \eqref{yli-size} and the fact from \eqref{abrush-card} that $\A_{brush}^{l_0}$ is non-empty, and \eqref{brush-bound} follows since
$\p(\A),\mu \gtrsim 1$ and $\lambda \lesssim 1$.

Now assume $\sigma \gg \delta^{-C/N} \delta$.  To prove \eqref{brush-bound} we will in fact prove the stronger bound
\be{x-strong}
 |\X'| \gtrapprox
\delta^{C/N} \mu \lambda^3 \sigma \delta^{n-2} \p(\A)^{-1}
\end{equation}
where
$$ \X' = \bigcup_{l \in \A_{brush}^{l_0}} Y_l^{(1)} \cap \Omega$$
and
$$ \Omega = \{ x \in \R^n: \delta^{C/N} \sigma \lesssim \dist(x,l_0) \lesssim \delta^{-C/N} \sigma \}.$$
From \eqref{yli-size}, \eqref{two-ends-three}, and elementary geometry we have
$$ |T_l \cap \X'| \approx \lambda \delta^{n-1}$$
for all $l \in \A_{brush}^{l_0}$.  Summing this in $l$ we obtain. 
$$
\sum_{l \in \A_{brush}^{l_0}} |T_l \cap \X'|
\approx |\A_{brush}^{l_0}| \lambda \delta^{n-1},$$
which we rewrite as
$$
\int_{\X'} \sum_{l \in \A_{brush}^{l_0}} \chi_{T_l \cap \Omega}
\approx |\A_{brush}^{l_0}| \lambda \delta^{n-1}.$$
We now use C\'ordoba's argument.  From Cauchy-Schwarz and the above we have
$$
|\X'|^{1/2} \|\sum_{l \in \A_{brush}^{l_0}} \chi_{T_l \cap \Omega}\|_2
\gtrapprox |\A_{brush}^{l_0}| \lambda \delta^{n-1}.$$
From this and \eqref{abrush-card}, it suffices to show that
\be{cord}
\|\sum_{l \in \A_{brush}^{l_0}} \chi_{T_l \cap \Omega}\|_2^2
\lessapprox \delta^{-C/N} |\A_{brush}^{l_0}| \delta^{n-1} \p(\A)^{-1},
\end{equation}
since \eqref{x-strong} then follows from algebra.

To prove \eqref{cord}, we expand the left-hand side as
$$ \sum_{l \in \A_{brush}^{l_0}} \sum_{l' \in \A_{brush}^{l_0}} |T_l \cap T_{l'} \cap \Omega|,$$
which we break up further as
$$ \sum_{\delta \lesssim \tau \lesssim 1}
\sum_{l \in \A_{brush}^{l_0}} \sum_{l' \in \A_{brush}^{l_0}: T_l \cap T_{l'} \cap \Omega \neq \emptyset,
\delta + |v(l)-v(l')| \sim \tau} |T_l \cap T_{l'} \cap \Omega|.$$
From elementary geometry we have
$$ |T_l \cap T_{l'}| \lesssim \delta^n \tau^{-1}.$$
It thus suffices to show that
$$ | \{ l' \in \A_{brush}^{l_0}: T_l \cap T_{l'} \cap \Omega \neq \emptyset, 
\delta + |v(l)-v(l')| \sim \tau \}| \lessapprox \delta^{-C/N} \p(\A) \tau/\delta$$
for each $l$, $\tau$.

Fix $l$, $\tau$.  The conditions $l' \in \A_{brush}^{l_0}$ and
$T_l \cap T_{l'} \cap \Omega \neq \emptyset$ force $l'$ to lie in a $\delta^{1-C/N}$-neighbourhood of the $2$-plane spanned by $l_0$ and (a slight translate of) $l$.  Together with the condition $\delta + |v(l)-v(l')|$, this constrains $T_{l'}$ to live in one of $O(\delta^{-C/N})$ boxes, each of dimension $C \times C\tau \times C\delta \times \ldots \times C\delta$.  The claim then follows from \eqref{plate-def}.
\end{proof}

\section{The bilinear reduction}

We now begin the proof of Theorem \ref{discrete-thm}.  

Fix $0 < \eps \ll 1$.  For each $0 < \delta \ll 1$, let $A(\delta) = A_\eps(\delta)$ denote the best constant such that
\be{linear} 
\| \sum_{l \in \A} \chi_{T_l} 
\|_{p'} \leq A(\delta) \delta^{-\frac{n}{p} + 1 - \eps}
m^{1/q - 1/r} (\delta^{n-1} |\A|)^{1/\qp}.
\end{equation}
for all choices of $m$, $\E$, $\E'$ and $\A$ satisfying \eqref{m-def}.
Clearly $A(\delta)$ is finite for each $\delta$; to prove Theorem \ref{discrete-thm}, we need to show
\be{a-delta}
A(\delta) \lessapprox 1.
\end{equation}

It will be convenient to denote the right-hand side of \eqref{linear} as
$\Q(\delta,\A)$, thus
\be{Q-def}
\Q(\delta,\A) = A(\delta) \delta^{-\frac{n}{p} + 1 - \eps}
m^{1/q - 1/r} (\delta^{n-1} |\A|)^{1/\qp}.
\end{equation}

By an inductive argument it suffices to prove \eqref{a-delta} assuming that
\be{local-max}
A(\delta) \sim \sup_{\delta \leq \delta' \ll 1} A(\delta')
\end{equation}

Fix $\delta$ so that \eqref{local-max} holds.  We may find $m$, $\E$, and
$\A$ such that
\be{linear-max} 
\| \sum_{l \in \A: v(l) \in \E} \chi_{T_l} 
\|_{p'} \sim \Q(\delta,\A).
\end{equation}

The estimate \eqref{linear-max} states that $\A$ is essentially an optimal configuration.  This has several consequences, at least heuristically.  Firstly, it implies that the generic angle between two lines in $\A$ is $\sim 1$.  Secondly, it implies a
``two-ends'' condition, which roughly asserts
that the contribution of the generic tube $T_l$ to \eqref{linear-max} is
not concentrated on a short interval.  We make these claims rigorous in the following sections, together with a technical uniformization reduction; these preliminaries will simplify the ensuing argument.  We remark that one needs $\eps > 0$ in order to obtain these reductions.

We begin with the assertion that the generic angle between two lines is $\sim 1$.  This is accomplished by

\begin{proposition} There exist subsets $\E_1, \E_2$ of $\E$ such that
\be{separation}
\dist(\E_1, \E_2) \sim 1
\end{equation}
and
\be{bilinear}
\| (\sum_{l \in \A: v(l) \in \E_1} \chi_{T_l}) (\sum_{l' \in \A: v(l') \in \E_2} \chi_{T_{l'}}) \|_{p'/2}^{1/2}
 \sim \Q(\delta,\A).
\end{equation}
\end{proposition}

Without \eqref{separation}, one could simply take $\E_1=\E_2=\E$ in the above
proposition.  The point of this proposition is that it allows one to restrict one's attention to pairs of tubes which intersect at large angle.  This bilinear reduction allows us to avoid many (but not all) of the difficulties
involving small angle intersections, which we have already encountered when
managing the $\sigma$ and $\tau$ parameters in the previous section.

\begin{proof}  By squaring \eqref{linear-max} we have
\be{square}
\| \sum_{l,l' \in \A} \chi_{T_l}  \chi_{T_{l'}}
\|_{p'/2} \sim \Q(\delta,\A)^2.
\end{equation}
Now let $0 < c_0 < 1$ be a small number to be chosen later, and consider the quantity
\be{narrow}
\| \sum_{l,l' \in \A: |v(l)-v(l')| < c_0} \chi_{T_l} \chi_{T_{l'}} \|_{p'/2}.
\end{equation}
Cover $\E$ by finitely overlapping sets $\E = \bigcup_\alpha \E_\alpha$ where each $\E_\alpha$ has diameter $O(c_0)$, and such that for every $v, v' \in \E$ with $|v-v'| \leq c_0$ there exists an $\alpha$ such that $v,v' \in \E_\alpha$.
We thus have
$$ \sum_{l,l' \in \A: |v(l)-v(l')| < c_0} \chi_{T_l} \chi_{T_{l'}}
\leq \sum_\alpha (\sum_{l \in \A_\alpha} \chi_{T_l})^2,
$$
where $\A_\alpha = \{ l \in \A: v(l) \in \E_\alpha\}$.
Since $p'/2 < 1$, we have the quasi-triangle inequality
\be{quasi}
 \| \sum_\alpha f_\alpha^2 \|_{p'/2} \leq (\sum_\alpha \| f_\alpha^2 \|_{p'/2}^{p'/2})^{2/p'} = (\sum_\alpha \|f_\alpha\|_{p'}^{p'})^{2/p'},
\end{equation}
(see e.g. \cite{TVV}), and so we may estimate \eqref{narrow} by
\be{narrow2}
(\sum_\alpha \| \sum_{l \in \A_\alpha} \chi_{T_l} \|_{p'}^{p'})^{2/p'}.
\end{equation}
We now claim that
\be{claim0}
\| \sum_{l \in \E_\alpha} \chi_{T_l} \|_{p'}
\lesssim c_0^{-(n-1)/p'} \Q(\delta/c_0, \A_\alpha)
\end{equation}
To see this, first apply a mild affine map to make $\E_\alpha$ centered at the origin, and apply the dilation $(\underline{x}, x_n) \to (\underline{x}/c_0, x_n)$, and then apply \eqref{linear} to the result; cf. \cite{TVV}.  

Since our choice of $p$, $q$ satisfy the scaling condition
$q = (n-1)p'$,
we may simplify \eqref{claim0} using \eqref{local-max} and \eqref{Q-def} to
$$
\| \sum_{l \in \A_\alpha} \chi_{T_l} \|_{p'}
\lesssim c_0^{\eps} (|\A_\alpha|/|\A|)^{1/\qp} \Q(\delta,\A).$$
Inserting this back into \eqref{narrow}
and using the elementary inequality
$$ \sum_\alpha (|\A_\alpha|/|\A|)^{p'/\qp} \leq (\sum_\alpha |\A_\alpha|/|\A|)^{p'/\qp}
\leq 1$$
which follows since $p' > \qp$, we obtain
$$ \eqref{narrow} \lesssim 
(c_0^{\eps} \Q(\delta,\A))^2.
$$
Comparing this with \eqref{square} we see that
$$
\| \sum_{l,l' \in \A: |v(l)-v(l')| \geq c_0} \chi_{T_l} \chi_{T_{l'}} \|_{p'/2}
\sim \Q(\delta,\A)^2$$
if we choose $c_0$ to be a sufficiently small number depending only on $n$
and $\eps$ (so $c_0 \sim 1$).

Now cover $\E$ by $O(c_0^{1-n})$ balls of diameter $c_0/4$.  By the pigeonhole principle and the above estimate we see that there must exist at least one pair $\E_1$, $\E_2$ of such balls with $\dist(\E_1,\E_2) \geq c_0/2$
such that
$$
\| \sum_{l,l' \in \A: v(l) \in \E_1, v(l') \in \E_2} \chi_{T_l} \chi_{T_{l'}} \|_{p'/2}
\gtrsim c_0^C
\Q(\delta,\A)^2.
$$
The claim follows.
\end{proof}

Note that the above argument is not restricted to this particular choice of
$p,q,r$.  See \cite{KLT}, \cite{TVV}, \cite{TV2} for variants of this argument.
The arguments in the next three sections are similarly not restricted to the
exponent choices in \eqref{squid}.

Henceforth $\E_1$, $\E_2$ will be fixed.

\section{Uniformity of multiplicity and density}

Let $\A$ be a subset of $\E \times \E'$ satisfying \eqref{m-def}, and let $\X$ be a subset of $\R^n$.
It would be convenient if we could ensure some uniformity on the multiplicity
function $\sum_{l \in A} \chi_{T_l}$ and the density function $|T_l \cap \X|$, as in Lemma \ref{tech}.  This is achieved by

\begin{lemma}\label{UNIF} 
Let $\A$ be a subset of $\E \times \E'$ satisfying \eqref{m-def}, and let $\X$ be a subset of $\R^n$.  Let $\lambda, \mu > 0$ be quantities satisfying \be{lambda-mu}
\mu |\X| = \lambda \delta^{n-1} |\A|.
\end{equation}
and
\be{norm-disc} \mu |\X|^{1/p'} \gtrapprox \Q(\delta,\A).
\end{equation}
Suppose $\X' \subset \X$, $\A' \subset \A$ are such that
\be{mu-hyp} \int_{\X'} \sum_{l \in \A'} \chi_{T_l} \approx \mu |\X|
\end{equation}
or equivalently that
\be{lambda-hyp} \sum_{l \in {\A'}} |T_l \cap {\X'}| \approx \lambda \delta^{n-1} |\A|.
\end{equation}
Then we have
\be{mu-u} 
\int_{x \in \X': \sum_{l \in \A'} \chi_{T_l}(x) \approx \mu}
\sum_{l \in \A'} \chi_{T_l}(x) \approx \mu |\X|
\end{equation}
and
\be{lambda-u}
 \sum_{l \in \A': |T_l \cap \X'| \approx \lambda \delta^{n-1}} 
|T_l \cap \X'| \approx \lambda \delta^{n-1} |\A|.
\end{equation}
Equivalently, we have
$$ | \{ x \in \X': \sum_{l \in \A'} \chi_{T_l}(x) \approx \mu \} | \approx |\X|$$
and
$$ | \{ l \in \A': |T_l \cap \X'| \approx \lambda \delta^{n-1} \} | \approx |\A|.$$
\end{lemma}

The condition \eqref{lambda-mu} is quite natural; cf. \eqref{mu-lambda-3}.  The condition \eqref{norm-disc} is a variant of \eqref{linear-max}, and states that $\mu |\X|^{1/p'}$ is essentially as large as possible.
Although this lemma is not phrased in a bilinear way, we will be able to combine it with the bilinear reduction (and the two-ends reduction in the next section) in Section \ref{unif-sec}.

\begin{proof}
We first prove \eqref{mu-u}.  Let $B = (\log (1/\delta))^\nu$, where $\nu$ is a large constant to be chosen later.
We trivially have
$$
\int_{x \in \X': \sum_{l \in \A'} \chi_{T_l}(x) \lessapprox B^{-1} \mu}
\sum_{l \in \A'} \chi_{T_l}(x) \lessapprox B^{-1} \mu |\X|.
$$
We now claim that
\be{b-high}
\int_{x \in \X': \sum_{l \in \A'} \chi_{T_l}(x) \gtrapprox B \mu}
\sum_{l \in \A'} \chi_{T_l}(x) \lessapprox B^{-(p'-1)} \mu |\X|;
\end{equation}
the claim then follows by subtracting these two estimates from \eqref{mu-hyp} and choosing $\nu$ suitably.

To prove \eqref{b-high}, we first observe that the left-hand side is bounded by
$$ \lessapprox (B\mu)^{1-p'}
\int (\sum_{l \in \A} \chi_{T_l})^{p'}.$$
By \eqref{linear} and \eqref{norm-disc}, this is bounded by
$$ \lessapprox (B\mu)^{1-p'} (\mu |\X|^{1/p'})^{p'},$$
and \eqref{b-high} follows.

Now we prove \eqref{lambda-u}, which is a dual of \eqref{mu-u}; the last two claims in the lemma then follow easily.  

As before we have
$$
 \sum_{l \in \A': |T_l \cap \X'| \lessapprox B^{-1} \lambda \delta^{n-1}} 
|T_l \cap \X'| \lessapprox B^{-1} \lambda \delta^{n-1} |\A|.
$$
It suffices to show that
\be{l-high}
 \sum_{l \in \A''} 
|T_l \cap \X'| \lessapprox {B'}^{-(q-1)} \lambda \delta^{n-1} |\A|.
\end{equation}
for all $B' \geq B$, where
$$ A'' = \{ l \in \A': |T_l \cap \X'| \approx B' \lambda \delta^{n-1} \};$$
by summing this for all dyadic $B' \geq B$ and using the exponential decay of the ${B'}^{-(q-1)}$ we can obtain the analogue of \eqref{b-high}.

Fix $B'$.  By definition of $A''$ we have
$$ \int_{\X'} \sum_{l \in \A''} \chi_{T_l} = \sum_{l \in \A''} |T_l \cap \X'|
\approx B' \lambda \delta^{n-1} |\A''|.$$
From H\"older we thus have
\be{fromhol} |\X|^{1/p} \| \sum_{l \in \A''} \chi_{T_l} \|_{p'}
\gtrapprox B' \lambda \delta^{n-1} |\A''|.
\end{equation}
From \eqref{linear} we have
$$ \| \sum_{l \in \A''} \chi_{T_l} \|_{p'} \leq \Q(\delta,\A'');$$
from \eqref{Q-def} and \eqref{norm-disc} we thus have
$$ \| \sum_{l \in \A''} \chi_{T_l} \|_{p'} \lessapprox
\mu |\X|^{1/p'} (|\A''|/|\A|)^{1/\qp}.$$
Inserting this into \eqref{fromhol} and using \eqref{lambda-mu} we obtain
$$ \lambda \delta^{n-1} |\A| (|\A''|/|\A|)^{1/\qp} \gtrapprox B \lambda \delta^{n-1} |\A''|,$$
which simplifies to
$$ |\A''| \lessapprox {B'}^{-q} |\A|,$$
and \eqref{l-high} follows from the definition of $\A''$.
\end{proof}

\section{The two ends reduction}

In order to apply Lemma \ref{three} we need (among other things) to obtain the conditions \eqref{lam-as} and \eqref{two-ends-three}.  The condition 
\eqref{lam-as} can essentially be guaranteed by Lemma \ref{UNIF}, but this lemma does not give us the two-ends condition \eqref{two-ends-three}.  To obtain this we shall use the following lemma.

\begin{lemma}\label{one-end}  Let $N \gg 1$, $\X$ be a subset of $\R^n$, and let $\A$ be a subset of $\E \times \E'$ satisfying \eqref{m-def}, and such that for every $l \in \A$ there exists an $x \in \R^n$ such that
$$ |T_l \cap \X \cap B(x,\delta^{1/N})| \gtrapprox \delta^{\eps/2N} |T_l \cap \X|.$$
Then we have
\be{conc-bound}
 \sum_{l \in \A} |T_l \cap \X| \lessapprox 
\delta^{\eps/2N} |\X|^{1/p} \Q(\delta,\A)
\end{equation}
\end{lemma}

The factor of $\delta^{\eps/N}$ in the above argument will allow us to conclude that for most tubes, the set $|T_l \cap \X|$ is not concentrated in a short end of the tube.  This type of ``two-ends condition'' first appears in \cite{W1}, \cite{W2}.

\begin{proof}
Cover $[0,1]$ by $\sim \delta^{-1/N}$ finitely overlapping intervals $I_\alpha$ of width
$\sim \delta^{1/N}$, and let $S_\alpha$ denote the slab $\R^{n-1} \times I_\alpha$.  For each $l \in \A$, we can then find an $\alpha = \alpha(l)$ such that
$$ |T_l \cap \X \cap S_\alpha| \gtrapprox \delta^{\eps/2N}
|T_l \cap \X|.$$
It thus suffices to show that
\be{conc2} 
\sum_\alpha \sum_{l \in \A_\alpha} |T_l \cap S_\alpha \cap \X|
\lessapprox \delta^{\eps/N} |\X|^{1/p} \Q(\delta,\A)
\end{equation}
where
$$ \A_{\alpha} = \{ l \in \A: \alpha(l) = \alpha \}.$$

Partition $\E$ into about $\delta^{(1-n)/N}$ refinements $\E_\beta$, each of which is $\delta^{1-1/N}$-separated.  We can split the left-hand side of \eqref{conc2} as
$$
\sum_\alpha
\sum_\beta
\int_{S_\alpha \cap \X} \sum_{l \in \A_{\alpha,\beta}} \chi_{T_l \cap \X}$$
where
$$ \A_{\alpha,\beta} = \{ l \in \A_\alpha: v(l) \in \E_\beta \}.$$
By H\"older, we may estimate this by
\be{asd}
\sum_\alpha
\sum_\beta |S_\alpha \cap \X|^{1/p}
\|  \sum_{l \in \A_{\alpha,\beta}} \chi_{T_l \cap \X} \|_{p'}.
\end{equation}
The sets $T_l \cap S_\alpha$ in the innermost sum can be rescaled to form a collection of $\delta^{1-1/N} \times 1$ tubes which continue to satisfy \eqref{m-def}.  Also, the set of directions $\E_\beta$ satisfies the correct separation condition for the scale $\delta^{1-1/N}$.
By a rescaled version of \eqref{linear} and \eqref{Q-def}, we can therefore bound the norm in \eqref{asd} by
$$
\lessapprox \delta^{n/Np'} \Q(\delta^{1-1/N},\A_{\alpha,\beta}),$$
which can be estimated using \eqref{Q-def}, \eqref{local-max} and algebra by
$$
\lessapprox \delta^{\eps/N} \Q(\delta,\A) \delta^{(n-1)/qN} (|\A_{\alpha,\beta}|/|\A|)^{1/\qp}.$$
Inserting this back into \eqref{asd}, we may estimate the left-hand side of
\eqref{conc2} as
$$ 
\lessapprox \delta^{\eps/N}
\Q(\delta,\A)
\sum_\alpha |S_\alpha \cap \X|^{1/p}
\delta^{(n-1)/qN} \sum_\beta (|\A_{\alpha,\beta}|/|\A|)^{1/\qp}.$$
Since we have $O(\delta^{(1-n)/N})$ $\beta$'s, we can use H\"older to obtain
$$
\delta^{(n-1)/qN} \sum_\beta (|\A_{\alpha,\beta}|/|\A|)^{1/\qp}
\lesssim (|\A_\alpha|/|\A|)^{1/\qp}.$$
We can thus bound the left-hand side of \eqref{conc2} as
$$ \lessapprox \delta^{\eps/N}
\Q(\delta,\A)
\sum_\alpha |S_\alpha \cap \X|^{1/p}
 (|\A_{\alpha}|/|\A|)^{1/\qp}.$$
By H\"older again, we bound this by
$$ \lessapprox \delta^{\eps/N}
\Q(\delta,\A)
(\sum_\alpha |S_\alpha \cap \X|^{q/p})^{1/q}$$
Since $q>p$, we can bound this by
$$ \lessapprox \delta^{\eps/N} \Q(\delta,\A) (\sum_\alpha |S_\alpha \cap \X|)^{1/p},$$
and \eqref{conc2} follows.
\end{proof}

\section{Plate number uniformization}\label{unif-sec}

We now combine the tools developed in the previous three sections to obtain the following technical uniformization lemma, which is analogous to Lemma \ref{tech}.  We use $\A_{i,0}$ for $i=1,2$ to denote the set
$$ \A_{i,0} = \{ l \in \A: v(l) \in \E_i\}.$$

\begin{lemma}\label{uniform}  Let the notation be as in the previous sections, and let $N \gg 1$ be a large number.  Then, if $\delta$ is sufficiently small depending on $\eps$ and $N$, there exist numbers $\mu, \lambda, \p_1, \p_2 > 0$
and sets
\be{anest} \A_i^{(3)} \subset \A_i^{(2)} \subset \A_i^{(1)} \subset \A_i^{(0)} = \A_{i,0}\hbox{ for } i=1,2,
\end{equation}
and
\be{xnest} \X^{(3)} \subset \X^{(2)} \subset \X^{(1)} \subset \X^{(0)} \subset \R^n
\end{equation}
such that 
\be{consistency}
|\X^{(0)}| \mu = |\A| \lambda \delta^{n-1}.
\end{equation}
and
\be{norm}
\mu |\X^{(0)}|^{1/p'} \approx
\Q(\delta,\A).
\end{equation}
Furthermore, one has
\be{lambda-unif}
|T_l \cap \X^{(j-1)}| \approx \lambda |T_l|
\end{equation}
\be{two-ends}
|T_l \cap \X^{(j-1)} \cap B(x,\delta^{1/N})| \lessapprox 
\delta^{\eps/2N} \lambda |T_l|
\end{equation}
for all $l \in \A_i^{(j)}$, $i = 1,2$, $j = 1,2,3$, $x \in \R^n$,
\be{mu-unif} \sum_{l \in \A_i^{(j)}} \chi_{T_l}(x) \approx \mu \hbox{ for all } x \in \X^{(j)},\ i = 1,2,\ j = 0,1,2,3,
\end{equation}
and
\be{p-unif}
\delta^{C/N} \p_i \lessapprox \p_i(\A_i^{(j)}) \lessapprox \delta^{-C/N} \p_i
\hbox{ for } i=1,2,\ j = 1,2,3.
\end{equation}
The implicit constants in these estimates may depend on $N$.
\end{lemma}

\begin{proof}  The first step is to find $\mu$ and $\X^{(0)}$.

Let $\mu_1, \mu_2$ range over all dyadic integers from $1$ to
$\delta^{-C}$.  Let $\X^{(0)}(\mu_1,\mu_2)$ denote the set
$$ \X^{(0)}(\mu_1,\mu_2) = \{ x: \sum_{l \in \A_{i,0}} \chi_{T_l}(x) \sim \mu_i \hbox{ for } i=1,2 \}.$$
Clearly we have
\be{bilinear-bound}
 \hbox{ LHS of \eqref{bilinear} } \sim
(\sum_{\mu_1} \sum_{\mu_2} \mu_1^{p'/2} \mu_2^{p'/2} |\X^{(0)}(\mu_1,\mu_2)|)^{1/p'}.
\end{equation}
Since the number of $\mu_1$ and $\mu_2$ is $\approx 1$, we can use the pigeonhole principle and conclude that there exist $\mu_1$, $\mu_2$ for which
\eqref{norm} holds with $\X^{(0)} = \X^{(0)}(\mu_1,\mu_2)$ and $\mu = (\mu_1 \mu_2)^{1/2}$. 

Fix this choice of $\mu_i$, $\mu$ and $\X^{(0)}$; this also fixes $\lambda$.  By construction we have
$$ \| \sum_{l \in \A_{i,0}} \chi_{T_l} \|_{p'} \gtrsim \mu_i |\X^{(0)}|^{1/p'}.$$
Combining this with \eqref{linear} we have
$$ \mu_i |\X^{(0)}|^{1/p'} \lesssim \Q(\delta,\A)$$
for $i=1,2$.
Combining this with \eqref{norm}
 we see that
$$ \mu_i \lessapprox \mu.$$
From the definition of $\mu$ we thus have $\mu_i \approx \mu$.  Since $\mu \lesssim \delta^{-C}$, we see from \eqref{norm}, \eqref{consistency} that
$|\X^{(0)}|, \lambda \gtrsim \delta^C$.

We now produce sets
$$ \X^{(0)} = \X_0 \supset \X_1 \supset \ldots \supset \X_{N^2}$$
and 
$$ \A_{i,0} \supset \A_{i,1} \supset \ldots \supset \A_{i,N^2}$$
with the properties that
\be{mass-bound}
|\X_k| \approx |\X_0|
\hbox{ for all } 0 \leq k \leq N^2,
\end{equation}
\be{lambda-iter}
|T_l \cap \X_{k-1}| \approx \lambda |T_l|,
\end{equation}
\be{two-ends-iter}
|T_l \cap \X_{k-1} \cap B(x,\delta^{1/N})| \lessapprox \delta^{\eps/2N} \lambda |T_l| 
\end{equation}
for all $l \in \A_{i,k}$, $i = 1,2$, $1 \leq k \leq N^2$, $x \in \R^n$,
and
\be{mu-unif-iter} \sum_{l \in \A_{i,k}} \chi_{T_l}(x) \approx \mu \hbox{ for all } x \in \X_k,\ i = 1,2,\ 0 \leq k \leq N^2.
\end{equation}

Clearly \eqref{mass-bound} and \eqref{mu-unif-iter} hold for $k=0$.
Now suppose inductively that $0 \leq k < N^2$ is such that
$\X_k, \A_{1,k}, \A_{2,k}$ have been constructed satisfying
\eqref{mu-unif-iter} and \eqref{mass-bound} for
this value of $k$.

We perform a certain sequence of dance steps.
From \eqref{mass-bound} and \eqref{mu-unif-iter} we have
$$\int_{X_k} \sum_{l \in \A_{1,k}} \chi_{T_l} \approx \mu |\X_0|,$$
which by \eqref{consistency} implies
$$
\sum_{l \in \A_{1,k}} |T_l \cap \X_k| \approx 
|\A| \lambda \delta^{n-1}.
$$
By Lemma \ref{UNIF} (noting that $\Q(\delta,\A_{1,k}) \leq \Q(\delta,\A)$; we shall need similar observations in the sequel), 
we thus have
\be{squish}
\sum_{l \in \A'_{1,k}} |T_l \cap \X_k| \approx 
|\A| \lambda \delta^{n-1}
\end{equation}
where $\A'_{1,k} \subseteq \A_{1,k}$ is the set
$$ \A'_{1,k} = \{ l \in \A_{1,k}: |T_l \cap \X_k| \approx \lambda \delta^{n-1} \}.$$

Now define the set $\A_{1,k+1} \subseteq \A'_{1,k}$ by
$$ \A_{1,k+1} = \{ l \in \A'_{1,k}:
|T_l \cap \X_k \cap B(x,\delta^{1/N})| \leq \delta^{\eps/2N} |T_l \cap \X_k|
\hbox{ for all } x \in \R^n. \}$$
From Lemma \ref{one-end} we have
$$
\sum_{l \in \A'_{1,k} \backslash \A_{1,k+1}} |T_l \cap \X_k| \lesssim 
\delta^{\eps/2N} |\X_0|^{1/p} \Q(\delta,\A);$$
by \eqref{norm} and \eqref{consistency} we thus have
$$
\sum_{l \in \A'_{1,k} \backslash \A_{1,k+1}} |T_l \cap \X_k| \lessapprox
\delta^{\eps/2N}
|\A| \lambda \delta^{n-1}.$$
Combining this with \eqref{squish} we obtain (if $\delta$ is sufficiently small)
$$
\sum_{l \in \A_{1,k+1}} |T_l \cap \X_k| \approx 
|\A| \lambda \delta^{n-1}.
$$
We may rewrite this using \eqref{consistency} as
$$
\int_{\X_k} \sum_{l \in \A_{1,k+1}} \chi_{T_l}(x) \approx
\mu |\X_0|.$$
By Lemma \ref{UNIF}, we have
$$ 
|\X'_k| \approx |\X_0|$$
where $\X'_k \subseteq \X_k$ is the set
$$ \X'_k = \{ x \in \X_k: \sum_{l \in \A_{1,k+1}} \chi_{T_l}(x)
\approx \mu\}.$$
In particular, from \eqref{mu-unif-iter} with $i=2$, we have
$$
\int_{\X'_k} \sum_{l \in \A_{2,k}} \chi_{T_l}(x) \approx
\mu |\X_0|.$$
By \eqref{consistency}, we may rewrite this as
$$
\sum_{l \in \A_{2,k}} |T_l \cap \X'_k| \approx 
|\A| \lambda \delta^{n-1}.$$
By Lemma \ref{UNIF} again, this implies
$$
\sum_{l \in \A'_{2,k}} |T_l \cap \X'_k| \approx |\A| \lambda \delta^{n-1}.
$$
where
$$ \A'_{2,k} = \{ l \in \A_{2,k}: |T_l \cap \X'_k| \approx \lambda \delta^{n-1} \}.$$
Defining
$$ \A_{2,k+1} = \{ l \in \A'_{1,k}:
|T_l \cap \X_k \cap B(x,\delta^{1/N})| \leq \delta^{\eps/2N} |T_l \cap \X_k|
\hbox{ for all } x \in \R^n. \}$$
we apply Lemma \ref{one-end}, \eqref{norm}, \eqref{consistency}, and the preceding estimate as before to conclude
$$
\sum_{l \in \A_{2,k+1}} |T_l \cap \X'_k| \approx |\A| \lambda \delta^{n-1}.$$
By \eqref{consistency} again, we rewrite this as
$$ \int_{\X'_k} \sum_{l \in \A_{2,k+1}} \chi_{T_l}(x) \approx \mu |\X_0|.$$
By Lemma \ref{UNIF} we have
$$ |\X_{k+1}| \approx |\X_0|$$
where
$$ \X_{k+1} = \{ x \in \X'_k: \sum_{l \in \A_{2,k+1}} \chi_{T_l}(x) \approx \mu \}$$
This completes the dance sequence.  One can easily verify that \eqref{mass-bound}, \eqref{two-ends-iter}, 
and \eqref{mu-unif-iter} are all satisfied for $k+1$ and $i=1,2$.  One now
replaces $k$ by $k+1$, and repeats the above dance.  Of course, the implicit constants in the bounds will depend on $k$ and hence on $N$.

The quantities $\p_i(\A_{i,k})$ are clearly monotone decreasing, and satisfy the trivial estimates $1 \lesssim \p_i(\A_{i,k}) \lesssim \delta^{-C}$.
By the pigeonhole principle one can then find $1 < k < N^2-1$ such that
$$ \p_i(\A_{i,k+2}) \geq \delta^{C/N} \p_i(\A_{i,k}) \hbox{ for } i=1,2.$$
The lemma then follows by setting $\X^{(j)} = \X_{k+j-1}$, $\A_i^{(j)} = \A_{i,k+j-1}$,
and $\p_i = \p_i(\A_{i,k})$ for $j=1,2,3$ and $i=1,2$.
\end{proof}

This argument can be extended to create arbitrarily longer sequences than the ones in the above lemma, but we shall not need to do so here.

\section{Estimates for a slab}

Let the notation be as in Lemma \ref{uniform}.  Define a $\theta$-slab to be a $\theta/2$-neighbourhood of a 2-plane in $\R^n$.

In the sequel we shall prove two propositions.

\begin{proposition}\label{slab}  Let $\delta \lesssim \theta \lesssim 1$, and let $S$ be a $\theta$-slab.  Then we have
\be{slab-est}
 |\X^{(1)} \cap S| \lessapprox \theta^{1/2} \lambda^{7/2 - n} 
|\A|^{\frac{n-2}{n-1}} m^{1/(n-1)} \delta^{n-2} \mu^{-1}.
\end{equation}
\end{proposition}

\begin{proposition}\label{hairbrush}  There exists a $\delta \lesssim \theta \lesssim 1$ and a $\theta$-slab $S$ such that
\be{hair-est} |\X^{(1)} \cap S| \gtrapprox \delta^{C/N} \mu \lambda^{7/2} m^{-1/2} \theta^{1/2} \delta^{n-2}.
\end{equation}
\end{proposition}

Suppose for the moment that both propositions were true.  Then we would have
$$ 
\delta^{C/N} \mu \lambda^{7/2} m^{-1/2} \delta^{n-2}
\lessapprox 
\lambda^{7/2 - n} |\A|^{\frac{n-2}{n-1}} m^{1/(n-1)} \delta^{n-2} \mu^{-1}.$$
If one uses \eqref{consistency} to eliminate $\lambda$, this becomes
(using \eqref{Q-def}, \eqref{squid} and a lot of algebra)
$$
\mu |\X^{(0)}|^{1/p'}
\lessapprox \delta^{C/N} \Q(\delta,\A)
$$
Comparing this with \eqref{norm} one obtains \eqref{a-delta} if $N$ is chosen sufficiently large depending on $\eps$.

It remains to prove the Propositions. 

\section{Proof of Proposition \ref{slab}}

We now prove Proposition \ref{slab}.  The estimate \eqref{slab-est} is not
best possible; it was chosen primarily so that it cancelled nicely against
\eqref{hair-est}.  Accordingly, our techniques shall be quite crude.  

Fix $\theta$ and $S$.  From \eqref{mu-unif} we have
$$ |\X^{(1)} \cap S| \approx \mu^{-1} \int_{\X^{(1)} \cap S} \sum_{l \in \A_1^{(0)}} \chi_{T_l}.$$
We can rewrite the right-hand side as
$$ \mu^{-1} \sum_{l \in \A_1^{(0)}} |\X^{(1)} \cap S \cap T_l|
\leq \mu^{-1} \sum_{l \in \A_1^{(0)}} |\X^{(0)} \cap S \cap T_l|.$$
For each $l$, let $\alpha(l)$ denote the quantity 
$$ \alpha(l) = \theta + \angle(l,S),$$
where $\angle(l,S)$ is the angle between $l$ and the plane in the middle of $S$.  From elementary geometry we have
$$ |S \cap T_l| \lesssim \delta^{n-1} \theta \alpha(l)^{-1},$$
and so by \eqref{lambda-unif} we have
$$ |\X^{(0)} \cap S \cap T_l| \lessapprox \delta^{n-1} \min(\theta \alpha(l)^{-1}, \lambda) \lessapprox \delta^{n-1} \theta^{1/2} \alpha(l)^{-1/2} \lambda^{1/2}.$$
Combining all these estimates we obtain
$$ |\X^{(1)} \cap S| \lessapprox
\mu^{-1} \delta^{n-1} \sum_{l \in \A_1^{(0)}} \theta^{1/2} 
\alpha(l)^{-1/2} \lambda^{1/2}.$$
From \eqref{lambda-unif} we have $\lambda \lessapprox 1$, so that
$\lambda^{1/2} \lessapprox \lambda^{7/2-n}$.  It thus suffices to show that
\be{slab-test}
\sum_{l \in \A_1^{(0)}} \alpha(l)^{-1/2}
\lessapprox \delta^{-1}
|\A|^{\frac{n-2}{n-1}} m^{\frac{1}{n-1}}.
\end{equation}
We can estimate the left-hand side of \eqref{slab-test} by
\be{chunk}
\sum_{\delta \lesssim \alpha \lesssim 1} \sum_{l \in \A: \alpha(l) \sim \alpha}
\alpha^{-1/2} \sim \sum_{\delta \lesssim \alpha \lesssim 1} \alpha^{-1/2}
|\{ l \in \A: \alpha(l) \sim \alpha\}|
\end{equation}
where $\alpha$ ranges over the dyadic numbers.  From \eqref{m-def} and the $\delta$-separated nature of $\E$ we have
$$ |\{ l \in \A: \alpha(l) \sim \alpha\}| \lesssim \alpha^{n-2} \delta^{1-n} m.$$
Interpolating this with the trivial bound of $|\A|$ we obtain
$$ |\{ l \in \A: \alpha(l) \sim \alpha\}| \lesssim \alpha^{\frac{n-2}{n-1}} \delta^{-1} m^{\frac{1}{n-1}} |\A|^{\frac{n-2}{n-1}}.$$
Inserting this back into \eqref{chunk} we obtain
\eqref{slab-test} since $\frac{n-2}{n-1} \geq \frac{1}{2}$.
This concludes the proof of Proposition \ref{slab-test}.
\endprf

It is clear that there is plenty of slack in the above estimate.  Indeed, the only time when \eqref{slab-est} is efficient is when $\lambda, \alpha, \theta \approx 1$, and when $|\E| \approx \delta^{1-n}$.  These phenomena seems to be a typical consequence of the two ends and bilinear reductions respectively.

\section{Proof of Proposition \ref{hairbrush}}

We now prove Proposition \ref{hairbrush}.  This shall be a modified version of
the hairbrush argument in \cite{W2}.

By symmetry we may assume
\be{p1-mono}
\p_1 \geq \p_2.
\end{equation}

Since $\p(\A_1^{(3)}) \gtrapprox \delta^{C/N} \p_1$ by \eqref{p-unif}, we see from \eqref{plate-def} that one can find a $\delta \lesssim w \lesssim 1$ and a $C \times Cw \times C\delta \times \ldots
\times C_\delta$ rectangle $R$ such that
\be{plate-card}
\frac{|\A_R|}{w/\delta} \gtrapprox \delta^{C/N} \p_1,
\end{equation}
where
$$ \A_R = \{ l \in \A_1^{(3)}: T_l \in R \}.$$
This rectangle $R$ shall form the stem of a hairbrush in $S \cap \X^1$.  Let $l_R$ denote the line generated by the first direction of $R$, and $\pi_R$ be the 2-plane generated by the first two directions of $R$; thus $R$ lies in the $C\delta$ neighbourhood of $\pi_R$ and in the $Cw$-neighbourhood of $l_R$.

By refining $\A_R$ slightly if necessary, we may assume that $w \ll \delta^{1/N}$; this may worsen the power of $\delta^{1/N}$ in \eqref{plate-card}, but is otherwise harmless.
From \eqref{separation} we thus have
\be{disj}
|v(l_R) - v(l)| \sim 1 \hbox{ for all } l \in \A_2^{(2)}.
\end{equation}

Since $\A_1^{(3)} \subset \E \times \E'$, we have from elementary geometry that
$$ |\A_R| \lesssim (w/\delta)^2.$$
Combining this with \eqref{plate-card} we see that
\be{w-bound}
w \gtrsim \delta^{C/N} \p_1 \delta.
\end{equation}

From \eqref{lambda-unif} we see that
\be{man}
 | T_l \cap \X^{(2)} | \approx \lambda \delta^{n-1}
\end{equation}
for all $l \in \A_R$.  From this we conclude

\begin{lemma}  We have
\be{x2-r}
 |\X^{(2)} \cap R| \gtrapprox \delta^{C/N}
\lambda^{3/2} w^{1/2} \p_1^{1/2} \delta^{n-3/2}.
\end{equation}
\end{lemma}

\begin{proof}
Firstly, from \eqref{plate-card} and elementary geometry we see that $\A_R$ must contain at least $\delta^{C/N} \p_1$ parallel lines, which with \eqref{man} and \eqref{p1-mono} gives
$$ |\X^{(2)} \cap R| \gtrapprox \delta^{C/N} \lambda \p_1 \delta^{n-1}.$$
It thus suffices to show
$$|\X^{(2)} \cap R| \gtrapprox \delta^{C/N} \lambda^2 |\A_R| \p_1^{-1} \delta^{n-1},$$
since \eqref{x2-r} follows by taking the geometric mean of these estimates and then using \eqref{plate-card}.

To prove this estimate we invoke C\'ordoba's argument as in the proof of \eqref{x-strong}.
Summing \eqref{man} over all $l \in \A_R$ we obtain
$$ \sum_{l \in \A_R} |T_l \cap \X^{(2)}| \approx \lambda \delta^{n-1} |\A_R|$$
which we rewrite as
$$ \int_{\X^{(2)} \cap R} \sum_{l \in \A_R} \chi_{T_l} \approx \lambda \delta^{n-1} |\A_R|.$$
By the Cauchy-Schwarz inequality we thus have
$$ |\X^{(2)} \cap R|^{1/2} \| \sum_{l \in \A_R} \chi_{T_l}\|_2 \gtrapprox
\lambda \delta^{n-1} |\A_R|.$$
It thus suffices to show that
\be{cord2} \| \sum_{l \in \A_R} \chi_{T_l}\|_2^2 \lessapprox \delta^{-C/N} \lambda^2 |\A_R| \p_1 \delta^{n-1}.
\end{equation}
Repeating the derivation of \eqref{cord}, we may estimate the left-hand side by
$$ \sum_{\delta \lesssim \tau \lesssim 1}
\sum_{l \in \A_R} \sum_{l' \in \A_R: T_l \cap T_{l'} \neq \emptyset, \delta + |v(l)-v(l')| \sim \tau} \delta^n \tau^{-1},$$
and the claim follows from the observation
$$ | \{ l' \in \A_R: T_l \cap T_{l'} \neq \emptyset, \delta + |v(l)-v(l')| \sim \tau \} | \lesssim \delta^{-C/N} \delta^{-1} \tau \p_1$$
which follows from \eqref{plate-def} and elementary geometry.
\end{proof}

Thus $\X^{(2)}$ has a large intersection with $R$.  We now wish to conclude that there are many tubes from $\A_2^{(2)}$ passing through $R$. 

Combining \eqref{x2-r} with \eqref{mu-unif} and \eqref{p1-mono} we have
$$ \int_R \sum_{l \in \A_2^{(2)}} \chi_{T_l}(x) \gtrapprox \lambda^{3/2} \mu w^{1/2} \p_2^{1/2} \delta^{n-3/2},$$
which we rewrite as
$$ \sum_{l \in \A_2^{(2)}} |T_l \cap R| \gtrapprox \lambda^{3/2} \mu w^{1/2} \p_2^{1/2} \delta^{n-3/2}.$$
For each dyadic $\delta \lesssim \theta \lesssim 1$,
let $\A_{brush}^\theta$ denote the set
$$ \A_{brush}^\theta = \{ l \in \A_2^{(2)}: T_l \cap R \neq \emptyset, \quad
\delta/w + \angle l, \pi_R \sim \theta \}.$$
We thus have
$$ \sum_{\delta/w \lesssim \theta \lesssim 1} \sum_{l \in \A_{brush}^\theta}
|T_l \cap R| \gtrapprox \lambda^{3/2} \mu w^{1/2} \p_2^{1/2} \delta^{n-3/2}.$$
By the pigeonhole principle, there must therefore exist a $\delta/w \lesssim \theta \lesssim 1$ such that
$$
\sum_{l \in \A_{brush}^\theta} |T_l \cap R| \gtrapprox \lambda^{3/2} \mu w^{1/2} \p_2^{1/2} \delta^{n-3/2}.$$
Fix this $\theta$.
From \eqref{disj} and the definition of $\A_{brush}^\theta$, we see from elementary geometry that
$$ | T_l \cap R| \lessapprox \delta^n \theta^{-1}.$$
Combining this with the previous, we see that
\be{abrush-est}
 |\A_{brush}^\theta| \gtrapprox \lambda^{3/2} \mu w^{1/2} \p_2^{1/2} \theta \delta^{-3/2}
\end{equation}
Thus to prove \eqref{hair-est} it suffices to show that
$$
|\X^{(1)} \cap S| \gtrapprox \delta^{C/N} \lambda^2 |\A_{brush}^\theta|
m^{-1/2} \theta^{-1/2} w^{-1/2} \p_2^{-1/2} \delta^{n-1/2}.
$$
We will in fact show the slightly stronger
\be{x1s}
|\X^{(1)} \cap S \cap \Omega| \gtrapprox \delta^{C/N} \lambda^{2} |\A_{brush}^\theta|
m^{-1/2} \theta^{-1/2} w^{-1/2} \p_2^{-1/2} \delta^{n-1/2}
\end{equation}
where $\Omega$ denotes the region
$$ \Omega = \{ x \in \R^n: \delta^{1/N} \lesssim \dist(x,l_R) \lesssim 1 \}.$$ 

We now foliate the hairbrush into three-dimensional regions in order to apply Lemma \ref{three}.

Let $S^{n-3}$ denote the portion of the unit sphere $S^{n-1}$ in $\R^n$ which is orthogonal to $\pi_R$, and let $\Gamma$ be a maximal 
$C^{-1} \delta$-separated subset of $S^{n-3}$.  For each $\omega \in \Xi$, let
$V_\omega$ denote the set
$$ V_\omega = \pi_R + \R \omega + B^n(0,C\delta);$$
these sets are $C\delta$-neighbourhoods of $3$-spaces.
From elementary geometry we may cover
$$ \A_{brush}^\theta = \bigcup_{\omega \in \Gamma} \A_{brush}^{\theta,\omega}$$
where
$$ \A_{brush}^{\theta,\omega} = \{ l \in \A_{brush}^\theta: T_l \subset V_\omega \}.$$
The sets $V_\omega \cap \Omega$ have an overlap of at most $O(\delta^{-C/N})$ as $\omega$ varies.  Thus
$$ |\X^{(1)} \cap S \cap \Omega| \gtrapprox \delta^{C/N}
\sum_{\omega \in \Gamma} |\X^{(1)} \cap S \cap V_\omega \cap \Omega|.$$
To show \eqref{x1s}, it thus suffices to show that
\be{x2s}
|\X^{(1)} \cap S \cap V_\omega \cap \Omega| \gtrapprox \delta^{C/N} \lambda^{2} |\A_{brush}^{\theta,\omega}|
m^{-1/2} \theta^{-1/2} w^{-1/2} \p_2^{-1/2} \delta^{n-1/2}
\end{equation}
for each $\omega \in \Gamma$.

Fix $\omega$.  The region $S \cap V_\omega \cap \Omega$ is essentially a $C \times C \times
C \theta \times C\delta \times \ldots \times C\delta$ box.  We cover this box
by about $w^{-1}$ smaller boxes $B_\alpha$ of dimensions $C \times C \times C w\theta \times C\delta \times \ldots \times C\delta$ such that $l_R$ is contained in the plane generated by the first two directions of this box.  Note that $w\theta \gtrsim \delta$ from the construction of $\theta$.  From
elementary geometry we see that for each $l \in \A_{brush}^{\theta,\omega}$
there exists a box $B_\alpha$ such that $T_l \subset B_\alpha$.  Also, the boxes
$B_\alpha$ have an overlap of $O(\delta^{-C/N})$.  Thus, by the same argument as before, it suffices to show that
\be{x3s}
|\X^{(1)} \cap B_\alpha| \gtrapprox \delta^{C/N} \lambda^{2} |\A_{brush}^{\theta,\omega,\alpha}|
m^{-1/2} \theta^{-1/2} w^{-1/2} \p_2^{-1/2} \delta^{n-1/2}
\end{equation}
where
$$ \A_{brush}^{\theta,\omega,\alpha} = \{ l \in \A_{brush}^{\theta,\omega}:
T_l \subset B_\alpha\}.$$

From \eqref{lambda-unif}, \eqref{two-ends} and elementary geometry we note that
$$ |T_l \cap \X^{(1)} \cap B_\alpha| \approx \lambda \delta^{n-1} \hbox{ for all } l \in \A_{brush}^{\theta,\omega,\alpha}.$$

Also, from elementary geometry we see that the set of directions
$\{ v(l): l \in \A_{brush}^{\theta,\omega,\alpha} \}$ is contained in a
$C \times Cw\theta \times C\delta \times \ldots \times C\delta$ box in
$B^{n-1}(0,1)$.  The claim \eqref{x3s} now follows from Lemma \ref{three},
and we are done.
\endprf

\end{document}